\documentclass[12pt]{amsart} 
 
\renewcommand{\hom}{{\rm Hom}}
\newcommand{\C}{\Bbb{C}}

\newcommand{\Z}{\Bbb{Z}}
\newcommand{\im}{{\rm im}}
\newcommand{\SL}{{\rm SL}}
\newcommand{\SU}{{\rm SU}}

\renewcommand{\H}{{\rm H}}
\newcommand{\HG}[1]{{\mathsf{H}}^1(#1,G)}

\newtheorem{theorem}{Theorem}

\title{Non-abelian Local Invariant Cycles}
\subjclass[2000]{14D05, 20F34, 55N20}

\author{Yen-lung Tsai}
\address{National Center for Theoretical Sciences, Mathematics Division,
 National Tsing Hua University, Hsinchu 300, Taiwan}
\email{yenlung@math.cts.nthu.edu.tw}
 
\author{Eugene Z. Xia}
\address{Department of Mathematics, National Cheng Kung University, Tainan 701, Taiwan}
\email{ezxia@ncku.edu.tw}

\thanks{Xia gratefully acknowledges partial support by National Science Council Taiwan grant NSC 93-2115-M-006-002.}

\begin{document}
 
\begin{abstract}
Let $f$ be a degeneration of K\"ahler manifolds. The local invariant cycle theorem states that for a smooth fiber of the degeneration, any cohomology class, invariant under the monodromy action, rises from a global cohomology class. Instead of the classical cohomology, one may consider the non-abelian cohomology. This note demonstrates that the analogous non-abelian version of the local invariant cycle theorem does not hold if the first non-abelian cohomology is the moduli space (universal categorical quotient) of the representations of the fundamental group.
\end{abstract}
 
\maketitle

A degeneration of K\"ahler manifolds is a proper map $f$ from a K\"ahler manifold $X$ onto the unit disk $\Delta$ such that $f$ is of maximum rank for all $s \in \Delta$ except at the point $s=0$. Let $\Delta^* = \Delta - \{ 0 \}$.  We call $X_t=f^{-1}(X_t)$ a smooth fiber or generic fiber when $t \in \Delta^*$ and $X_0=f^{-1}(0)$ the singular or degenerated fiber. We assume the singularity in $X_0$ is of normal crossing. 

Fix $t \in \Delta^*$ and a base point $x \in X_t$ once and for all.  There is a monodromy action (see, for example, \cite{kulikov:hodge})
\[
\pi_1(\Delta^*) \times \H^n(X_t,\C) \to \H^n(X_t,\C).
\]
The local invariant cycle theorem states that a cohomology class in $\H^n(X_t,\C)$, fixed by the monodromy action, is a restriction of a cohomology class in $\H^n(X,\C)$ \cite{clemens:sequence, schmid:vhs, steenbrink:limit}. Fix a generator $T \in \pi_1(\Delta^*) \cong \Z$. Then $T$ determines the monodromy action and gives rise to the isomorphisms
\[
T^*: \H^n(X_t,\C) \to \H^n(X_t,\C),   \ \ \ T_* : \pi_1(X_t, x) \to \pi_1(X_t,
y).
\]
These isomorphisms are actually induced from a Picard-Lefschetz diffeomorphism which we shall also denote by
\[
T: X_t \to X_t,
\]
where $y = T(x)$.

The map $f$ induces a strong deformation retract $X \to X_0$, hence, also an isomorphism $\H^n(X_0, \C) \to \H^n(X, \C)$ \cite{clemens:sequence, kulikov:hodge, persson:degenerations}.  Define 
\[
c : X_t \hookrightarrow X \to X_0,
\]
where the first map is the inclusion and the second the strong deformation retract.  Then $c$ induces homomorphisms:
\[
c^*: \H^n(X_0, \C) \to \H^n(X_t,\C), \ \ \ c_* : \pi_1(X_t,x) \to \pi_1(X_0,c(x)).
\] 
The local invariant cycle theorem states that \cite{clemens:sequence, schmid:vhs, steenbrink:limit}
\begin{equation} \label{csa}
\im (c^*) = \{ \alpha \in \H^n(X_t,\C) \mid T^*(\alpha) = \alpha \}.
\end{equation} 

To generalize the above statement to the non-abelian context, we choose the universal categorical quotient 
\[
\HG{X_s} = \hom(\pi_1(X_s,x),G)/\!/G, \ \ s \in \Delta
\]
to be the first non-abelian cohomology of $X_s$, where $G$ is an algebraic group acting on 
$\hom(\pi_1(X_s,x),G)$ by conjugation \cite{lubotzky:representation, mumford:GIT}. For any $\rho \in \hom(\pi_1(X_t,x),G)$, 
denote by $[\rho]$ its image in $\HG{X_t}$ by the universal morphism.

The map $T_*$ induces a morphism
\[
T^\# : \hom (\pi_1(X_t,T(x)),G) \to \hom (\pi_1(X_t,x),G)
\]
defined by
\[
T^\#(\rho)(A) = \rho \circ T_*^{-1} (A)
\]
for all $A \in \pi_1(X_t,x)$. Similarly, $c$ induces a map
\[
c^\# : \hom(\pi_1(X_0,c(x)),G) \to \hom(\pi_1(X_t,x),G)
\]
defined by
\[
c^\#(\rho)(A) = \rho \circ c_*(A).
\]
The maps $T^\#$ and $c^\#$ descend to morphisms on their respective
universal categorical quotients:
\[
T^\# : \HG{X_t} \to \HG{X_t}, \ \ \ c^\# : \HG{X_0} \to \HG{X_t}.
\]
In this setting, the analogous statement of the local invariant cycle theorem is
\begin{equation} \label{cs}
\im (c^\#) = \{ [\rho] \in \HG{X_t} \mid T^\#([\rho]) =
[\rho] \}.
\end{equation}
Notice that when $G=\C$ (the addition group), $\HG{X_s}$ is the regular
first cohomology $\H^1(X_s, \C)$ of $X_s$ for $s \in \Delta$ and (\ref{cs}) reduces to the classical local
invariant cycle theorem (\ref{csa}). 

\begin{theorem}\label{main}
For each $g>1$, there exists an $n$ and a degeneration $f: X \to \Delta$ with the generic fiber $X_t$
a Riemann surface of genus $g$ and such that 
\[
\im (c^\#) \neq \{ [\rho] \in \HG{X_t} \mid T^\#([\rho]) =
[\rho] \},
\]
where $G = \SL(n,\C)$.  More precisely, for such an $n$, there exists an irreducible representation
$\rho \in
\hom(\pi_1(X_t,x),G)$ with $\im(\rho) \subset \SU(n)$ 
and $[\rho] \in \{ [\rho] \in \HG{X_t} \mid T^\#([\rho]) = [\rho] \}]$ but $[\rho] \not\in \im (c^\#)$.
\end{theorem}
\begin{proof}
For each $n$, let $G = \SL(n,\C)$ and $K = \SU(n)$.
Let $\Gamma$ be the mapping class group of $X_t$.  Then each
$\gamma \in \Gamma$ induces a map
\[
\gamma : \hom(\pi_1(X_t,\gamma(x)),G) \to \hom(\pi_1(X_t,x),G) : \ \
\gamma(\rho)(A) = \rho \circ \gamma^{-1}(A).
\]
These maps provide a $\Gamma$-action on $\HG{X_t}$.
In \cite{andersen:fixedPoint}, Andersen proved that for $g > 1$, there
exists infinitely many $n$ such that there exists irreducible $\rho \in
\hom(\pi_1(X_t),G)$ with $\im(\rho) \subset K$ and $\gamma([\rho]) =
[\rho]$ for all $\gamma \in \Gamma$.  

Now suppose our choices of $\rho$ and $n$
satisfy the hypothesis and conclusions of Andersen's theorem \cite{andersen:fixedPoint}.  Then
since $T : X_t \to X_t$ is a
diffeomorphism, $T_* \in \Gamma$.  Hence $T^\#([\rho]) = [\rho]$.

The fundamental group $\pi_1(X_t,x)$ has a presentation
\[
\langle A_i,  1 \le i \le 2g, \mid \prod_{i = 1}^g A_i A_{g+i}
A_i^{-1} A_{g+i}^{-1}\rangle,
\]
where the $A_i$'s correspond to simple loops beginning and ending at
$x$.  Since $\rho$ is an irreducible representation, $\rho(A_i) \neq I \in
G$ for some $i$.

Now let $f : X \to \Delta$ be a degeneration with $A_i$ as its
vanishing cycle.  Then $c_*(A_i) = e \in \pi(X_0, c(x))$.  Hence 
for every $[\alpha] \in \HG{X_0}$, 
\[
c^\#(\alpha)(A_i) = \alpha(c_*(A_i)) = I \in G.
\]
Since $\rho(A_i) \neq I$, 
$[\rho] \neq [\alpha]$ for all $[\alpha] \in \im(c^\#)$.
\end{proof}

We refer to \cite{andersen:fixedPoint} for the detailed and explicit descriptions of these irreducible
representation classes fixed by the $\Gamma$-action.  Here we simply mention that all such representations have finite images in $\SU(n)$. It will be interesting to see whether there exist examples
$[\rho]$, satisfying the conclusion of Theorem~\ref{main}, with $\im(\rho)$ 
being Zariski dense in $G$.  Unfortunately, such representations are unlikely to satisfy the overly stringent condition of being fixed by the $\Gamma$-action. In fact, it is likely that the
$\Gamma$-orbit of $[\rho]$ is Zariski dense if $\im(\rho)$ is Zariski dense in $G$.  This is the case when $G = \SL(2,\C), K = \SU(2)$ \cite{previte:topologicalDyanmics1, previte:topologicalDynamics2}.

\end{document}